\documentclass[reqno]{amsart}
\title[Identities involving partitions with distinct even parts]{Identities involving partitions with distinct even parts and $4$-regular partitions}
\usepackage{amssymb,amsmath,amsthm,epsfig,graphics,latexsym,hyperref}
\theoremstyle{definition}
\newtheorem{definition}{Definition}
\theoremstyle{plain}

\newtheorem{theorem}    {Theorem}

\newtheorem{corollary}  {Corollary}

\theoremstyle{remark}

\numberwithin{equation}{section}

\newcommand{\fr}{\frac}

 \newmuskip\pFqskip
\mathchardef\pFcomma=\mathcode`, 

 \newmuskip\pGqskip
\mathchardef\pGcomma=\mathcode`, 

\begin{document}
\author[ G. E. Andrews and M. El Bachraoui]{George E. Andrews and Mohamed El Bachraoui}
\address{The Pennsylvania State University, University Park, Pennsylvania 16802}
\email{andrews@math.psu.edu}
\thanks{First author partially supported by Simons Foundation Grant 633284}
\address{Dept. Math. Sci,
United Arab Emirates University, PO Box 15551, Al-Ain, UAE}
\email{melbachraoui@uaeu.ac.ae}
\keywords{integer partitions, overpartition pairs, $q$-series, Bailey lemma.}
\subjclass[2000]{11P81; 05A17; 11D09}
\begin{abstract}
It is well known that the number of partitions into distinct even parts equals the number of $4$-regular partitions.
In this paper we prove identities relating certain restricted partitions into distinct even parts
with restricted $4$-regular partitions.
\end{abstract}
\date{\textit{\today}}
\maketitle
\section{Introduction}\label{sec introduction}
Throughout let $q$ be a complex number satisfying $|q|<1$, let $m$ and $n$ be nonnegative integers.
We will use the following standard notation for $q$-series~\cite{Andrews, Gasper-Rahman}
\[
(a;q)_0 = 1,\  (a;q)_n = \prod_{j=0}^{n-1} (1-aq^j),\quad
(a;q)_{\infty} = \prod_{j=0}^{\infty} (1-aq^j),
\]
\[
(a_1,\ldots,a_k;q)_n = \prod_{j=1}^k (a_j;q)_n,\ \text{and\ }
(a_1,\ldots,a_k;q)_{\infty} = \prod_{j=1}^k (a_j;q)_{\infty}.
\]
We will need the following basic facts of $q$-series~\cite{Andrews, Gasper-Rahman}
\begin{equation}\label{basic-facts}
(a;q)_{n+m} = (a;q)_{m} (aq^{m};q)_n,\
(a;q)_{\infty} = (a;q)_n (aq^n;q)_{\infty},\
(a;q)_{\infty} = (a;q^2)_{\infty}(aq;q^2)_{\infty}
\end{equation}
along with the $q$-binomial theorem
\begin{equation}\label{q-binomial}
\sum_{n=0}^\infty \fr{(a;q)_n}{(q;q)_n} z^n = \fr{(az;q)_\infty}{(z;q)_\infty}.
\end{equation}
%
%
Recall that for any nonnegative integer $l$, by $l$-regular partitions we mean the partitions whose summands are not divisible by $l$.
Let ${\rm ped}(n)$ denote the number of partitions of $n$ with distinct even parts~\cite{Andrews 2009}.
These partitions and their arithmetic properties have been studied extensively in recent years, see for
instance~\cite{And-Hir-Sel, Chen, Merca}.
One clearly has~\cite{Andrews 2009}
\[
\sum_{n=0}^\infty {\rm ped}(n) q^n = \fr{(-q^2;q^2)_\infty}{(q;q^2)_\infty} = \fr{(q^4;q^4)}{(q;q)_\infty}
\]
from which it follows that ${\rm ped}(n)$ equals the number of partitions of $n$ whose parts are not divisible by $4$.
See~\cite{Honsberger} for some other similar results.
Our main goal in this paper is to prove more identities relating restricted partitions with distinct even parts and
restricted $4$-regular partitions.

The paper is organized as follows.
In Section~\ref{sec results} we introduce our partitions and state the main results.
In Sections~\ref{sec proof thm main-1}-\ref{sec proof thm main-3} we give the proofs of the main results and
Section~\ref{sec conclusion} is devoted to some remarks suggested by this work.
\section{Definitions and main results}\label{sec results}
\begin{definition}\label{DE1}
Let ${\rm DE1}(n)$ denote the number of partitions of $n$ in which no even part is repeated and the largest part is odd.
Then it easy to see that
\[
\sum_{n=0}^\infty{\rm DE1}(n) q^n = \sum_{n=0}^\infty \fr{(-q^2;q^2)_n q^{2n+1}}{(q;q^2)_{n+1}}.
\]
\end{definition}
We shall prove the following two identities.
\begin{theorem}\label{thm main-1}
There holds
\begin{align}
\sum_{n\geq 0} q^{2n} (q^{4n+4};q^4)_\infty (q;q)_{2n} &= \fr{2(q^4;q^4)_\infty}{1+q}- \fr{(q;q)_\infty}{1+q} \label{help id-1}. \\
(1+q) \sum_{n\geq 0} \fr{(-q^2;q^2)_n q^{2n+1}}{(q;q^2)_{n+1}} &= \fr{(q^4;q^4)_\infty}{(q;q)_\infty}-1 \label{main id-1}.
\end{align}
\end{theorem}
We get the following immediate consequence of~\eqref{main id-1}.
\begin{corollary}\label{cor main-1}
For $n>0$, ${\rm DE1}(n)+{\rm DE1}(n-1)$ equals the number of the $4$-regular partitions of $n$.
\end{corollary}
For example, we have ${\rm DE1}(8)=9$ counting
\[
7+1, 5+3, 5+2+1, 5+1+1+1, 3+3+2, 3+3+1+1, 3+2+1+1+1, 3+1+1+1+1+1,
\]
\[
1+1+1+1+1+1+1+1
\]
and ${\rm DE1}(7)=7$ counting
\[
7, 5+2, 5+1+1, 3+3+1, 3+2+1+1, 3+1+1+1+1, 1+1+1+1+1+1+1.
\]
Furthermore, there are $16$ partitions for $n=8$ which are $4$-regular. Namely,
\[
7+1, 6+2, 6+1+1, 5+3, 5+2+1, 5+1+1+1, 3+3+2, 3+3+1+1, 3+2+2+1,
\]
\[
3+2+1+1+1, 3+1+1+1+1+1, 2+2+2+2,
2+2+2+1+1, 2+2+1+1+1+1,
\]
\[
2+1+1+1+1, 1+1+1+1+1+1+1+1.
\]
We now introduce our second  partitions.
\begin{definition}\label{DE2}
Let ${\rm DE2}(n)$ denote the number of partitions of $n$ in which no even part is repeated, the largest part is odd and appears at least twice.
Then it easy to see that
\[
\sum_{n=0}^\infty{\rm DE2}(n) q^n = \sum_{n=0}^\infty \fr{(-q^2;q^2)_n q^{4n+2}}{(q;q^2)_{n+1}}.
\]
\end{definition}
We will prove the following two identities.
\begin{theorem}\label{thm main-2}
There holds
\begin{align}
\sum_{n\geq 0} q^{2n} (q^{4n+4};q^4)_\infty (q;q)_{2n+1} &= \fr{2(1-q)(q^4;q^4)_\infty}{1+q^3}- \fr{(q;q)_\infty}{1+q^3} \label{help id-2}. \\
(1+q^3) \sum_{n\geq 0} \fr{(-q^2;q^2)_n q^{4n+2}}{(q;q^2)_{n+1}} &= \fr{(q^4;q^4)_\infty}{(q^2;q)_\infty}-1 \label{main-id-2}.
\end{align}
\end{theorem}
As the right hand-side of~\eqref{main-id-2} generates the number of partitions of $n$ into parts each $>1$ and not divisible by $4$,
we get the following immediate corollary of Theorem~\ref{thm main-2}.
\begin{corollary}\label{cor main-2}
For $n>0$, ${\rm DE2}(n)+{\rm DE2}(n-3)$ equals the number of $4$-regular partitions of $n$ into parts each $>1$.
\end{corollary}
For example, for $n=10$ we have ${\rm DE2}(10)=5$ counting
\[
5+5, 3+3+3+3+1, 3+3+2+1+1, 3+3+1+1+1+1, 1+1+1+1+1+1+1+1+1+1,
\]
${\rm DE2}(7)=2$ counting
\[
3+3+1, 1+1+1+1+1+1+1
\]
and there are $7$ partition of $10$ in which each part $>1$ and not divisble by $4$:
\[
10, 7+3, 6+2+2, 5+5, 5+3+2, 3+3+2+2, 2+2+2+2+2.
\]
We now present our third  partition function.
\begin{definition}\label{DE3}
Let ${\rm DE3}(n)$ denote the number of partitions of $n$ in which no even part is repeated, the largest part is odd and appears exactly once.
Then it easy to see that
\[
\sum_{n=0}^\infty{\rm DE3}(n) q^n = \sum_{n=0}^\infty \fr{(-q^2;q^2)_n q^{2n+1}}{(q;q^2)_{n}}.
\]
\end{definition}
\begin{theorem}\label{thm main-3}
There holds
\begin{align}
\sum_{n\geq 0} q^{4n+1} (q^{4n+4};q^4)_\infty (q;q)_{2n}
&=\fr{2q^2(q^4;q^4)_\infty}{1+q^3}+ \fr{q(1-q) (q;q)_\infty}{1+q^3} \label{help id-3}. \\
(1+q^3) \sum_{n\geq 0} \fr{(-q^2;q^2)_n q^{2n+1}}{(q;q^2)_{n}} &= \fr{q^2(q^4;q^4)_\infty}{(q;q)_\infty}-q^2+q \label{main-id-3}.
\end{align}
\end{theorem}
We have the following consequence of~\eqref{main-id-3}.
\begin{corollary}\label{cor main-3}
For $n>1$, ${\rm DE3}(n+2)+{\rm DE3}(n-1)$ equals the number of $4$-regular partitions of $n$.
\end{corollary}
For example, ${\rm DE3}(10)=11$ enumerating
\[
9+1, 7+3, 7+2+1, 7+1+1+1, 5+4+1, 5+3+2, 5+3+1+1, 5+2+1+1+1,
\]
\[
5+1+1+1+1+1, 3+2+1+1+1+1+1, 3+1+1+1+1+1+1+1
\]
and ${\rm DE3}(7)=5$ enumerating.
\[
7, 5+2, 5+1+1, 3+2+1+1, 3+1+1+1+1
\]
Furthermore, we have already seen above that the number of $4$-regular partitions of $8$ is $16$.
We close this section with the following consequence of Corollary~\ref{cor main-1} and Corollary~\ref{cor main-3}.
\begin{corollary}
For $n>1$ we have
\[
{\rm DE3}(n+2)+{\rm DE3}(n-1) = {\rm DE1}(n)+{\rm DE1}(n-1).
\]
\end{corollary}
\section{Proof of Theorem~\ref{thm main-1}}\label{sec proof thm main-1}
We will apply the following formula of Andrews, Subbarao, and Vidyasagar~\cite{And-Sub-Vid}
\begin{equation}\label{Andr-Subb-Vidy}
\sum_{n\geq 0} \fr{(a;q)_n}{(b;q)_n} q^n
=\fr{q(a;q)_\infty}{b(b;q)_\infty \big(1-\fr{aq}{b}\big)} + \fr{1-\fr{q}{b}}{1-\fr{aq}{b}}.
\end{equation}
Letting in~\eqref{Andr-Subb-Vidy} $q\to q^2$, $a=q$, and $b=-q^2$,
we deduce
\[
\sum_{n\geq 0} q^{2n} (q^{4n+4};q^4)_\infty (q;q)_{2n}
= (q^4;q^4)_\infty \sum_{n\geq 0} \fr{(q;q)_{2n} q^{2n}}{(q^4;q^4)_n}
\]
\[
=(q^4;q^4)_\infty \sum_{n \geq 0} q^{2n} \fr{(q;q^2)_n}{(-q^2;q^2)_n}
\]
\[
=(q^4;q^4)_\infty \Big( \fr{-(q;q^2)_\infty}{(1+q)(-q^2;q^2)_\infty} + \fr{2}{1+q} \Big)
\]
\[
=\fr{2(q^4;q^4)_\infty}{1+q}- \fr{(q;q^2)_\infty (q^2;q^2)_\infty}{1+q}
\]
\[
=\fr{2(q^4;q^4)_\infty}{1+q}- \fr{(q;q)_\infty}{1+q},
\]
which proves~\eqref{help id-1}. We now prove~\eqref{main id-1}.
With the help of~\eqref{q-binomial}, we get
\[
\sum_{n\geq 0} \fr{(q;q)_{2n} q^{2n}}{(q^4;q^4)_n}
= (q;q)_\infty \sum_{n\geq 0} \fr{q^{2n}}{(q^4;q^4)_n (q^{2n+1};q)_{\infty}}
\]
\[
=(q;q)_\infty \sum_{n\geq 0} \fr{q^{2n}}{(q^4;q^4)_n} \sum_{m\geq 0} \fr{q^{2nm + m}}{(q;q)_m}
\]
\[
=(q;q)_\infty \sum_{m\geq 0} \fr{q^{m}}{(q;q)_m}  \sum_{n\geq 0} \fr{q^{(2m+2)n}}{(q^4;q^4)_n}
\]
\[
=(q;q)_\infty \sum_{m\geq 0} \fr{q^{m}}{(q;q)_m}  \fr{1}{(q^{2m+2};q^4)_\infty}
\]
\[
=(q;q)_\infty \sum_{m\geq 0} \Big(\fr{q^{2m}}{(q;q)_{2m} (q^{4m+2};q^4)_\infty} +  \fr{q^{2m+1}}{(q;q)_{2m+1} (q^{4m+4};q^4)_\infty} \Big)
\]
\[
=\fr{(q;q)_\infty}{(q^2;q^4)_\infty}  \sum_{m\geq 0} \fr{q^{2m}(q^2;q^4)_m}{(q;q)_{2m}}
+ \fr{(q;q)_\infty}{(q^4;q^4)_\infty}  \sum_{m\geq 0} \fr{q^{2m+1}(q^4;q^4)_m}{(q;q)_{2m+1}}
\]
\[
=\fr{(q;q)_\infty}{(q^2;q^4)_\infty}  \sum_{m\geq 0} \fr{q^{2m}(-q;q^2)_m}{(q^2;q^2)_{m}}
+ \fr{(q;q)_\infty}{(q^4;q^4)_\infty}  \sum_{m\geq 0} \fr{q^{2m+1}(-q^2;q^2)_m}{(q;q^2)_{m+1}}
\]
\[
=\fr{(q;q)_\infty}{(q^2;q^4)_\infty}\fr{(-q^3;q^2)_\infty}{(q^2;q^2)_\infty}
+ \fr{(q;q)_\infty}{(q^4;q^4)_\infty}  \sum_{m\geq 0} \fr{q^{2m+1}(-q^2;q^2)_m}{(q;q^2)_{m+1}}
\]
\[
=\fr{1}{1+q} +  \fr{(q;q)_\infty}{(q^4;q^4)_\infty}  \sum_{m\geq 0} \fr{q^{2m+1}(-q^2;q^2)_m}{(q;q^2)_{m+1}}.
\]
Now multiply both sides of the foregoing identity by $\fr{(q^4;q^4)_\infty}{(q;q)_\infty}$ to obtain
\[
 \fr{(q^4;q^4)_\infty}{(q;q)_\infty} \sum_{n\geq 0}\fr{(q;q)_{2n} q^{2n}}{(q^4;q^4)_n}
=  \fr{(q^4;q^4)_\infty}{(1+q)(q;q)_\infty} +  \sum_{m\geq 0} \fr{q^{2m+1}(-q^2;q^2)_m}{(q;q^2)_{m+1}}
\]
or equivalently,
\[
\sum_{m\geq 0} \fr{q^{2m+1}(-q^2;q^2)_m}{(q;q^2)_{m+1}}
=\fr{1}{(q;q)_\infty}\sum_{n\geq 0} q^{2n} (q;q)_{2n} (q^{4n+4};q^4)_\infty -  \fr{(q^4;q^4)_\infty}{(1+q)(q;q)_\infty}.
\]
Finally combining the previous identity with~\eqref{help id-1}, we arrive at
\[
\sum_{m\geq 0} \fr{q^{2m+1}(-q^2;q^2)_m}{(q;q^2)_{m+1}}
=\fr{2(q^4;q^4)_\infty}{(1+q)(q;q)_\infty}- \fr{1}{1+q} - \fr{(q^4;q^4)_\infty}{(1+q)(q;q)_\infty}
\]
\[
=\fr{(q^4;q^4)_\infty}{(1+q) (q;q)_\infty}- \fr{1}{1+q},
\]
which gives the desired formula after multiplying both sides by $1+q$.
\section{Proof of Theorem~\ref{thm main-2}}\label{sec proof thm main-2}
Noting that
\[
\sum_{n\geq 0} q^{2n} (q^{4n+4};q^4)_\infty (q;q)_{2n+1}
=  (q^4;q^4)_\infty \sum_{n=0}^\infty \fr{(q;q)_{2n+1}}{(q^4;q^4)_n} q^{2n}
\]
\[
=(1-q)(q^4;q^4)_\infty \sum_{n=0}^\infty \fr{(q^3;q^2)_n}{(-q^2;q^2)_n} q^{2n},
\]
formula~\eqref{help id-2} follows by an application of~\eqref{Andr-Subb-Vidy}
to $q\to q^2$, $a=q^3$, and $b=-q^2$.
Furthermore, with the help of~\eqref{q-binomial} we find
\[
\sum_{n\geq 0} \fr{(q;q)_{2n+1} q^{2n}}{(q^4;q^4)_n}
= (q;q)_\infty \sum_{n\geq 0} \fr{q^{2n}}{(q^4;q^4)_n (q^{2n+2};q)_{\infty}}
\]
\[
=(q;q)_\infty \sum_{n\geq 0} \fr{q^{2n}}{(q^4;q^4)_n} \sum_{m\geq 0} \fr{q^{2nm + 2m}}{(q;q)_m}
\]
\[
=(q;q)_\infty \sum_{m\geq 0} \fr{q^{2m}}{(q;q)_m}  \sum_{n\geq 0} \fr{q^{(2m+2)n}}{(q^4;q^4)_n}
\]
\[
=(q;q)_\infty \sum_{m\geq 0} \fr{q^{2m}}{(q;q)_m}  \fr{1}{(q^{2m+2};q^4)_\infty}
\]
\[
=(q;q)_\infty \sum_{m\geq 0} \Big(\fr{q^{4m}}{(q;q)_{2m} (q^{4m+2};q^4)_\infty} +  \fr{q^{4m+2}}{(q;q)_{2m+1} (q^{4m+4};q^4)_\infty} \Big)
\]
\[
=\fr{(q;q)_\infty}{(q^2;q^4)_\infty}  \sum_{m\geq 0} \fr{q^{4m}(q^2;q^4)_m}{(q;q)_{2m}}
+ \fr{(q;q)_\infty}{(q^4;q^4)_\infty}  \sum_{m\geq 0} \fr{q^{4m+2}(q^4;q^4)_m}{(q;q)_{2m+1}}
\]
\[
=\fr{(q;q)_\infty}{(q^2;q^4)_\infty}  \sum_{m\geq 0} \fr{q^{4m}(-q;q^2)_m}{(q^2;q^2)_{m}}
+ \fr{(q;q)_\infty}{(q^4;q^4)_\infty}  \sum_{m\geq 0} \fr{q^{4m+2}(-q^2;q^2)_m}{(q;q^2)_{m+1}}
\]
\[
=\fr{(q;q)_\infty}{(q^2;q^4)_\infty}\fr{(-q^5;q^2)_\infty}{(q^4;q^2)_\infty}
+ \fr{(q;q)_\infty}{(q^4;q^4)_\infty}  \sum_{m\geq 0} \fr{q^{4m+2}(-q^2;q^2)_m}{(q;q^2)_{m+1}}
\]
\[
=\fr{1-q}{1+q^3} +  \fr{(q;q)_\infty}{(q^4;q^4)_\infty}  \sum_{m\geq 0} \fr{q^{4m+2}(-q^2;q^2)_m}{(q;q^2)_{m+1}}.
\]
Now multiply both sides of the foregoing identity by $\fr{(q^4;q^4)_\infty}{(q;q)_\infty}$ to obtain
\[
\fr{(q^4;q^4)_\infty}{(q;q)_\infty} \sum_{n\geq 0}\fr{(q;q)_{2n+1} q^{2n}}{(q^4;q^4)_n}
=  \fr{(1-q)(q^4;q^4)_\infty}{(1+q^3)(q;q)_\infty} +  \sum_{m\geq 0} \fr{q^{4m+2}(-q^2;q^2)_m}{(q;q^2)_{m+1}}
\]
which combined with~\eqref{help id-2} yields
\[
\sum_{m\geq 0} \fr{q^{4m+2}(-q^2;q^2)_m}{(q;q^2)_{m+1}}
=\fr{2(q^4;q^4)_\infty}{(1+q^3)(q^2;q)_\infty}- \fr{1}{1+q^3} - \fr{(q^4;q^4)_\infty}{(1+q^3)(q^2;q)_\infty}
\]
\[
=\fr{(q^4;q^4)_\infty}{(1+q^3)(q^2;q)_\infty}- \fr{1}{1+q^3}.
\]
Finally multiply both sides by $1+q^3$ to obtain the desired formula.
\section{Proof of Theorem~\ref{thm main-3}}\label{sec proof thm main-3}
By~\eqref{help id-1} and~\eqref{help id-2}
\[
\sum_{n\geq 0} q^{4n+1} (q^{4n+4};q^4)_\infty (q;q)_{2n}
=\sum_{n\geq 0} q^{2n} (q^{4n+4};q^4)_\infty (q;q)_{2n} - \sum_{n\geq 0} q^{2n} (q^{4n+4};q^4)_\infty (q;q)_{2n+1}
\]
\[
=\fr{2q^2(q^4;q^4)_\infty}{1+q^3}+ \fr{q(1-q) (q;q)_\infty}{1+q^3},
\]
which confirms~\eqref{help id-3}.
In addition, we find with the help of~\eqref{q-binomial}
\[
\sum_{n\geq 0} \fr{(q;q)_{2n} q^{4n}}{(q^4;q^4)_n} = (q;q)_\infty \sum_{n\geq 0}\fr{q^{4n}}{(q^4;q^4)_n (q^{2n+1};q)_\infty}
\]
\[
= (q;q)_\infty \sum_{n\geq 0}\fr{q^{4n}}{(q^4;q^4)_n} \sum_{m\geq 0}\fr{q^{2nm+m}}{(q;q)_m}
\]
\[
=(q;q)_\infty \sum_{m\geq 0} \fr{q^m}{(q;q)_m (q^{2m+4};q^4)_\infty}
\]
\[
=(q;q)_\infty \sum_{m\geq 0} \fr{q^{2m}}{(q;q)_{2m} (q^{4m+4};q^4)_\infty}
+(q;q)_\infty \sum_{m\geq 0} \fr{q^{2m+1}}{(q;q)_{2m+1} (q^{4m+6};q^4)_\infty}
\]
\[
=\fr{(q;q)_\infty}{(q^4;q^4)_\infty}\sum_{m\geq 0} \fr{(q^4;q^4)_m q^{2m}}{(q;q)_{2m}}
+\fr{(q;q)_\infty}{(q^6;q^4)_\infty}\sum_{m\geq 0} \fr{(q^6;q^4)_m q^{2m+1}}{(q;q)_{2m+1}}
\]
\[
=\fr{(q;q)_\infty}{(q^4;q^4)_\infty}\sum_{m\geq 0} \fr{(-q^2;q^2)_m q^{2m}}{(q;q^2)_{m}}
+\fr{q(q;q)_\infty}{(1-q)(q^6;q^4)_\infty}\sum_{m\geq 0} \fr{(-q^3;q^2)_m q^{2m}}{(q^2;q^2)_{m}}
\]
\[
=\fr{(q;q)_\infty}{(q^4;q^4)_\infty}\sum_{m\geq 0} \fr{(-q^2;q^2)_m q^{2m}}{(q;q^2)_{m}} + \fr{q}{1+q^3}.
\]
Therefore multiplying both sides by $\fr{(q^4;q^4)_\infty}{(q;q)_\infty}$ gives
\[
\fr{(q^4;q^4)_\infty}{(q;q)_\infty} \sum_{n\geq 0} \fr{(q;q)_{2n} q^{4n}}{(q^4;q^4)_n}
=\fr{q(q^4;q^4)_\infty}{(1+q^3)(q;q)_\infty} + \sum_{m\geq 0} \fr{(-q^2;q^2)_m q^{2m}}{(q;q^2)_{m}}
\]
which by~\eqref{help id-3} means that
\[
\fr{2q(q^4;q^4)_\infty}{(1+q^3) (q;q)_\infty} + \fr{1-q}{1+q^3}
= \fr{q(q^4;q^4)_\infty}{(1+q^3) (q;q)_\infty} + \sum_{m\geq 0} \fr{(-q^2;q^2)_m q^{2m}}{(q;q^2)_{m}}.
\]
Hence, rewriting and multiplying the foregoing formula we arrive at
\[
(1+q^3) \sum_{m\geq 0} \fr{(-q^2;q^2)_m q^{2m}}{(q;q^2)_{m}}
=\fr{q(q^4;q^4)_\infty}{(q;q)_\infty} + 1-q.
\]
This completes the proof.
\section{Concluding remarks}\label{sec conclusion}
Our proofs for Theorems~\ref{thm main-1}-\ref{thm main-3} and their corollaries rely on the theory of basic hypergeometric series.
Therefore, it is natural to ask for combinatorial proofs for these results.

\bigskip
\noindent{\bf Acknowledgment.} The authors are grateful to the referee for valuable comments and interesting suggestions. 

\noindent{\bf Data Availability Statement.\ }
Not applicable.
\end{document}